%% file: main.tex
\documentclass[letterpaper, 10 pt, conference]{ieeeconf}  

\IEEEoverridecommandlockouts                              
\usepackage{graphics} 
\usepackage{epsfig} 
\usepackage{subfigure}
\usepackage{dblfloatfix}
\usepackage{times} 
\usepackage{amsmath} 
\usepackage{amssymb}  
\usepackage{xcolor}

\title{\LARGE \bf
Congestion-aware Bi-modal Delivery Systems Utilizing Drones
}
\author{Mark Beliaev$^{1}$, Negar Mehr$^{2}$, and Ramtin Pedarsani$^{3}$
\thanks{$^{1}$Mark Beliaev is with the Department of Electrical and Computer Engineering,
        University of California Santa Barbara, Santa Barbara, CA 93106, USA
        {\tt\small mbeliaev@ucsb.edu}}%
\thanks{$^{2}$Negar Z. Mehr is with the Department of Aerospace Engineering,
        University of Illinois at Urbana-Champaign, Champaign, IL 61820, USA
        {\tt\small negar@illinois.edu}}%
\thanks{$^{3}$Ramtin Pedarsani is with the Department of Electrical and Computer Engineering,
        University of California Santa Barbara, Santa Barbara, CA 93106, USA
        {\tt\small ramtin@ucsb.edu}}%
}
\input{macros.tex}

\begin{document} 

\maketitle
\thispagestyle{empty}
\pagestyle{empty}
\begin{abstract}
Bi-modal delivery systems are a promising solution to the challenges posed by the increasing demand of e-commerce. Due to the potential benefit drones can have on logistics networks such as delivery systems, some countries have taken steps towards integrating drones into their airspace. In this paper we aim to quantify this potential by developing a mathematical model for a Bi-modal delivery system composed of trucks and drones. We propose an optimization formulation that can be efficiently solved in order to design socially-optimal routing and allocation policies. We incorporate both societal cost in terms of road congestion and parcel delivery latency in our formulation. Our model is able to quantify the effect drones have on mitigating road congestion, and can solve for the path routing needed to minimize the chosen objective. To accurately capture the effect of stopping trucks on road latency, we model it using SUMO by simulating roads shared by trucks and cars. Based on this, we show that the proposed framework is computationally feasible to scale due to its reliance on convex quadratic optimization techniques. 
\end{abstract}
\section{Introduction}
The growing demand of e-commerce is causing a precipitous increase in delivery truck traffic, impacting travel time for both public and private vehicles within densely populated urban areas \cite{pavone_24}. For the past decade, e-commerce sales have grown at a yearly rate of $15\%$, with an astonishing $40\%$ increase in 2020 alone \cite{e_commerce}. It has been shown that trucks negatively impact road congestion by changing the behaviour of passenger car drivers around them \cite{impact_truck1}. Some have looked at ways to address these impacts, finding that the best strategy is restricting all heavy vehicles from using the road \cite{impact_truck2}. In addition to increased road congestion, the upsurge in truck traffic has led to other negative consequences such as air pollution and parking space shortages \cite{air_pollution,parking}. To address these issues, researchers have been developing Bi-modal delivery systems that incorporate parcel-delivering drones \cite{pavone, pavone_29}.\par
\begin{figure}[t]
\vspace{5pt}
  \centering
  \includegraphics[width=3.2in]{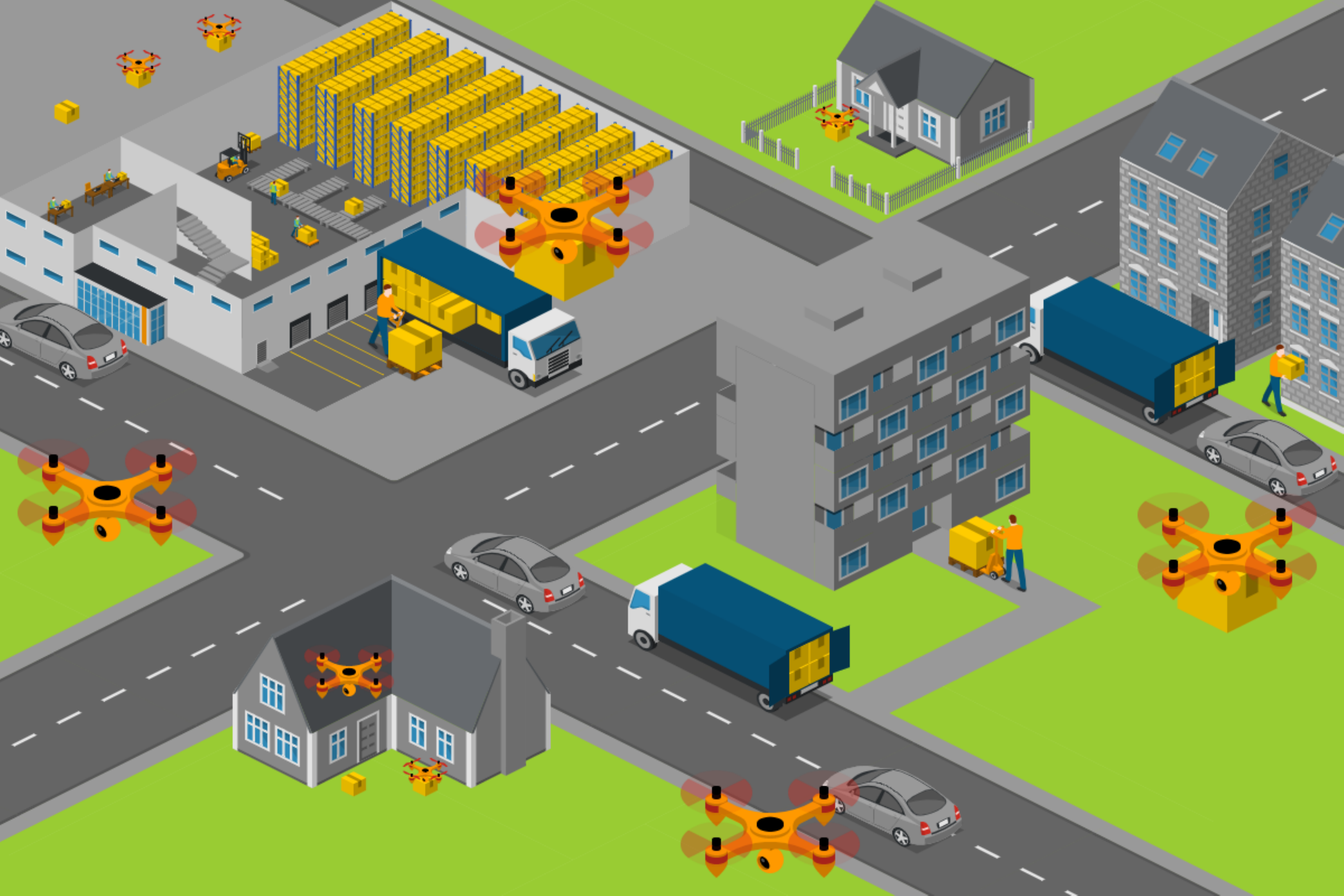}
  \caption{A depiction of a Bi-modal delivery system utilizing both trucks and drones. Delivery trucks making stops along roads with less lanes are shown to have a greater effect on road congestion. On the other hand, drones are shown landing directly at their assigned delivery destinations, hence not contributing negative effects on road congestion.}
  \label{fig:front1}
  \vspace{-20pt}
\end{figure}

Due to the aerial reach of Unmanned Aerial Vehicles, or drones, they possess an inherent versatility that makes them promising for logistics networks \cite{pavone_28}. Europe has already taken critical steps towards integrating drones into their current airspace infrastructure \cite{uspace}. In their effort to overcome the challenges posed to drone delivery services, Parcel at Deutsche Post DHL Group has demonstrated success in multiple field tests, using their recent Parcelcopter 4.0 to deliver medicine for the 400,000 residents of the Ukerewe island district of Lake Victoria \cite{dhl_drones}. These examples show the feasibility of implementing drones within logistics networks such as parcel delivery systems, depicted in Fig. \ref{fig:front1}.\par
In this paper we aim to quantify the potential drones have in alleviating traffic congestion induced by parcel delivery services by developing a mathematical model for a Bi-modal delivery system composed of trucks and drones. We set out to propose an optimization formulation that can be efficiently solved in order to design socially-optimal routing and allocation policies. Our goal is to incorporate both societal cost in terms of road congestion and parcel delivery latency in our optimization formulation.\par
We propose a new network model consisting of two separate graphs which represent the road network and the aerial network. A requisite for our model is the ability to capture the average time it takes parcels to be delivered. In addition, we want our model to capture the societal impact that different routing and allocation policies have on the road network. This way, given an objective function posed as a trade-off between societal cost and parcel latency, we form an optimization problem to solve for the socially-optimal routing strategy.\par
To model the effect of stopping delivery trucks on an individual link within the road network we use the open source traffic simulation package “SUMO” \cite{SUMO2018}. Based on the simulation results, we can model the latency of a link as linear with respect to the flow of stopping and non-stopping vehicles. This way we can develop a scalable optimization framework by making our objective function quadratic and convex.\par
The main contributions of the paper are as follows: 
\begin{itemize}
    \item We propose a new Bi-modal delivery network model consisting of delivery trucks, drones, and regular vehicles. 
    \item We describe a scalable convex optimization framework which solves for the socially-optimal allocation and routing scheme for our network model. 
    \item  We study the effect stopping trucks have on roads using SUMO simulations, and derive a mathematical latency function for roads shared between regular vehicles and delivery trucks.
    \item We perform numerical evaluations on our results to show the effect drone delivery systems have on road congestion and parcel delivery time.
\end{itemize}

\section{Related Work}
In the following, we discuss a few areas of work related to ours.
\subsection{Drone Delivery Systems} 
The application of drones within delivery systems has received notable interest \cite{drone_overview}. Lee \cite{drone_modular} has looked at the limitations posed to drones by their battery capacity and limited payload. Other work has looked at the facility location problem, optimizing for the locations of drone refueling stations in a large urban area \cite{drone_refuel}. More recently, a case study of Amazon prime air in the city of San Francisco has demonstrated how a proposed heuristic algorithm can approximately solve this NP-hard facility location problem \cite{drone_amazon}. Although these works consider the management and design of drone delivery systems as a whole, the optimizations performed do not consider the actual routing of drones within the network and their impact on the traffic network.\par
\subsection{Vehicle Routing Problem}
Solving for optimal routing schemes of a parcel delivery system within a road network can be viewed as an example of a vehicle routing problem (VRP) \cite{pavone_8,pavone_37_uav_route}. These are typically solved with mixed integer linear programming (MILP) formulations which scale poorly in real world applications. On top of this, Dorling et al. \cite{drone_vrp} have shown that naively applying existing VRPs to systems utilizing drones can lead to inefficient results. One notable work formulates the hybrid vehicle-drone routing problem for pick-up and delivery services as a MILP \cite{rev1_karak2019hybrid}, but they do not consider the affect the proposed routes have on society, and they do not allow trucks to dispatch and deliver packages. Although there are other works that consider relevant variants of the VRP \cite{rev3_liu2019hybrid,rev7_vidal2020concise}, the goal in these problems is to minimize the time it takes to complete the required routes, with no concern for the impact these schedules have on society.\par
\subsection{Multi-modal Route Planning}
There has been considerable progress within the past decade in developing new techniques for route planning in transportation networks \cite{route_plan_survey, rev4_bergmann2020integrating,rev5_56goswami2019one,rev8_reis2013rail}. Specifically, we have seen work that incorporates drones in the traffic model and optimizes for their routes, while allowing drones to utilize multiple modes of transportation \cite{pavone_10,pavone,drone_cost,rev2_ham2018integrated}. Other works have considered using drones in tandem with delivery trucks, without exploiting public transit \cite{pavone_18,pavone_36}. Although related to our problem, none of these works consider the effect drones have on the transportation network from a societal perspective. By decoupling the routes drones can take within the aerial network from the routes delivery trucks can take within the road network, we focus on the impact drones can have on road congestion as well as parcel latency within the road network.\par
\begin{figure}[b]
  \centering
  \vspace{-20pt}
  \includegraphics[width=2.8in]{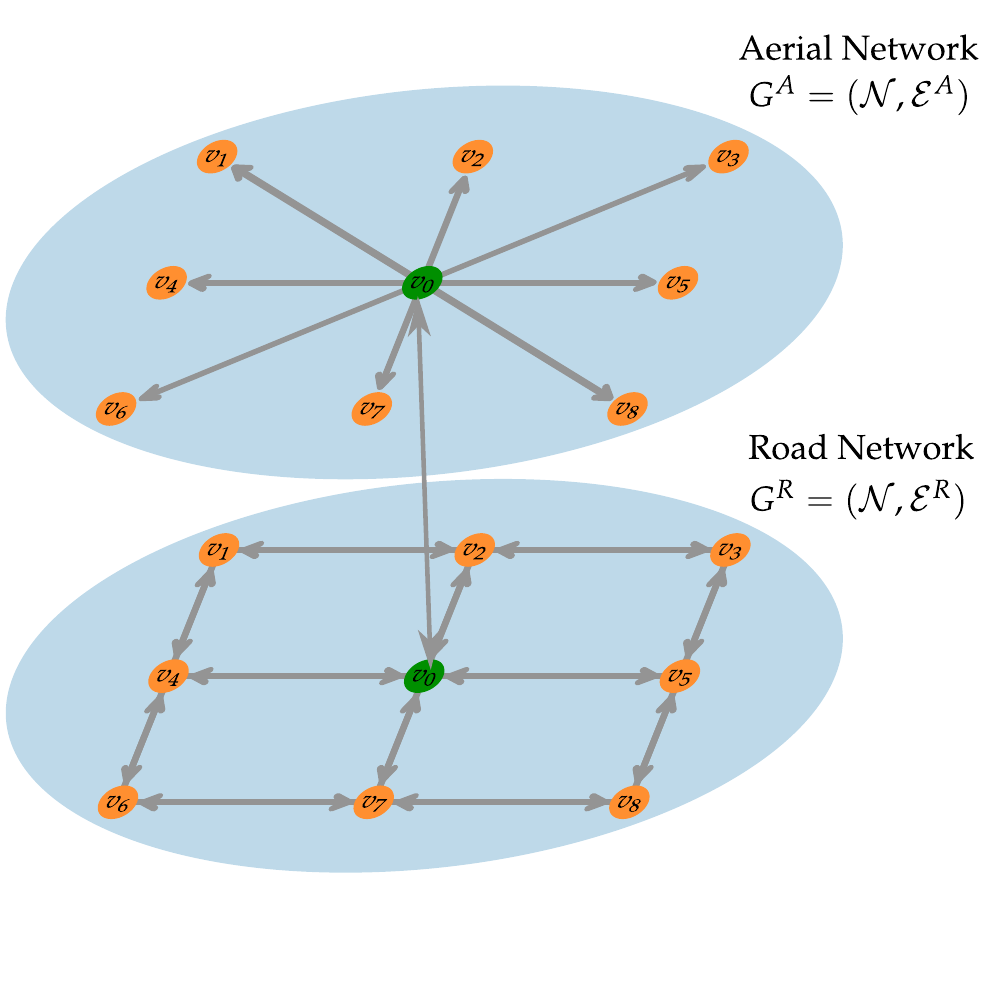}
  \caption{This figure visualizes the proposed network model. We see that the general structure of the directed graph $\Graph^R$ allows us to model typical transportation networks. On the other hand, the graph $\Graph^A$ is given the topology of a star graph since we allow drones to take direct paths to their destinations within the aerial network.}
  \label{fig:network}
\end{figure}
\section{Network Model}
In our work, we model the delivery modes as two separate networks, without interlinked paths between them. Each network is modeled by a directed graph composed of the same set of nodes $\Nodes$ and two different sets of edges $\Links^R$ and $\Links^A$, where $\Graph^R=(\Nodes,\Links^R)$ and $\Graph^A=(\Nodes,\Links^A)$ represent the road and aerial networks respectively. We assume that each edge $\link\in\Links^R,\Links^A$ connects two distinct nodes $(v,w)$, representing a directed edge from $v$ towards $w$. Additionally, we assume that both graphs are connected and denote node $\cnode$ as the main delivery hub, making all other nodes possible delivery destinations. Finally, we restrict each edge $\link\in\Links^A$ to only using $\cnode$ as the origin node, since the versatility of drones should allow them to have direct links between the main delivery hub and the destinations. We depict these two graphs in Fig. \ref{fig:network} for clarity.\par

To simplify notation, we introduce the sets of paths $\Paths^R$ and $\Paths^A$, where $\p\in\Paths^R$ denotes a path within $\Graph^R$ and $\p\in\Paths^A$ denotes an aerial path within $\Graph^A$. For the road network, each path $\p\in\Paths^R$ connects the delivery hub $\cnode$ to one of the remaining nodes using a collection of edges $\link\in\Links^R$. For the aerial network, each path $\p\in\Paths^A$ only contains one edge $\link\in\Links^A$ since $\Graph^A$ is a star network. Furthermore, for each node $\node\in\Nodes$, we define the subsets $\uplink_\node^R,\downlink_\node^R\subset\Links^R$ which correspond to the set of edges arriving at $\node$ and the set of edges departing from $\node$, respectively.\par

\subsection{Demand}
In our model, we let $\demand_\node$ represent the parcel demand of node $\node\in\Nodes$, measured in parcels per hour. Note that we actually use the set $\Nodes\setminus\{\cnode\}$ since the delivery hub does not have any demand, but we leave this distinction to the reader in order to avoid cumbersome notation. To ensure package delivery at all destinations, for each node $\node\in\Nodes$, we must have:
\begin{equation}
    \demand_\node = \demand_\node^R + \demand_\node^A,
    \label{eq:vertex_demand}
\end{equation}
where $\demand_\node^R$ and $\demand_\node^A$ correspond to the components of node $\node$'s  parcel demand satisfied by trucks within the road network and drones within the aerial network, respectively. How these components are derived is explained in the following section.\par

\subsection{Flow}
In contrast to the model's delivery demands measured in parcels per hour, we express the individual path and edge flows in units of vehicles per hour. For the road network, we use $\flow_\p^T$ to represent the flow of parcel-carrying trucks along a given path $\p\in\Paths^R$. Moreover, we let $\flow^T$ be the vector of truck flows along all the paths $p\in\Paths^R$. For edge $\link\in\Links^R$, we use $\flow_\link^T=\sum_{\p\in\Paths^R:\link\in\p}\flow_\p^T$ to denote the flow of trucks along edge $\link$. Additionally, each edge $\link\in\Links^R$ has a corresponding flow of non-delivering vehicles $\flow_\link^0$ which are constants defined by the model. We refer to $\flow_\link^0$ as the \textit{nominal} vehicular flow along edge $\link$ to outline that it is not a decision variable.\par 
For the road network, we assume that each truck carries $\ppertruck$ parcels. Thus, at each node $\node$, the number of parcels delivered by trucks is: 
\begin{equation}
    \tfrac{1}{\ppertruck}\demand_\node^R = \sum_{\link\in\uplink_\node^R}\flow_\link^T - \sum_{\link\in\downlink_\node^R}\flow_\link^T. 
    \label{eq:flow_cons}
\end{equation}
We can see from \eqref{eq:flow_cons} that $\demand_\node^R$ is measured as the \textit{decrease} in truck flow passing through node $\node$, since the difference between incoming and outgoing delivery trucks determines the number of parcels dropped off at node $\node$. Note that we do not consider the effect returning trucks have on the traffic network as they no longer need to make stops and are typically a negligible portion of the overall flow on any given edge.\par 

For the aerial network, we let $\flow_\p^D$ denote the flow of drones along path $\p\in\Paths^A$. Due to the limited payload of drones \cite{pavone_14}, we assume that each drone delivers only one package. This assumption coupled with the aforementioned restriction which allows only one aerial path to each destination node, makes flow $\flow_\p^D$ equivalent in magnitude to parcel demand $\demand_\node^A$ when $\p$ is the aerial path that arrives at destination node $\node$. We note that $\flow^T$ is used as our decision variable, deriving the value for $\flow_\p^D$ from \eqref{eq:vertex_demand} and \eqref{eq:flow_cons} as the leftover parcel demand that needs to be delivered through drones.\par

\subsection{Parcel Latency}
We would like to derive the mean latency experienced by parcels in the system. We express this average latency in units of hours as:
\begin{equation}
    \Latency = \frac{\Latency^R + \Latency^A}{\sum\limits_{\node\in\Nodes}\demand_\node},
    \label{eq:avg_latency}
\end{equation}
where we let $\Latency^R$ and $\Latency^A$ represent the sums of experienced latencies for parcels delivered within the road network and aerial network respectively. This way, dividing by the total demand in \eqref{eq:avg_latency} gives us the average latency we want. We go on to describe how these two components of latency are derived.\par

To derive the latency component $\Latency^R$, we first need to compute the individual edge latencies $\latency_\link^R$ for all edges $\link\in\Links^R$. Once these are computed, we can express the sum of latencies experienced by parcels delivered via trucks as:
\begin{equation}
    \Latency^R = \sum_{\p\in\Paths^R}\ppertruck\flow_\p^T\Big(\sum_{\link\in\p}\latency_\link^R\Big),
    \label{eq:truck_latency}
\end{equation}
where the term in the parenthesis is the latency along a path $\p$. As will be shown later, the edge latencies $\latency_\link^R$ are a linear function of the two flow components on that edge, making the expression for $\Latency^R$ quadratic with respect to our control variable $\flow^T$. Details of how this function was calculated are saved for Section \ref{sec:SUMO} describing our SUMO simulations. Note that \eqref{eq:truck_latency} is only valid when all trucks are carrying the same number of parcels $\ppertruck$. As long as the demand of individual vertices is large enough, the majority of trucks will be full, making this inaccuracy negligible.\par
For the drone network, we note that drone flows $\flow_\p^D$ do not have an effect on the latencies of the aerial paths since in the regimes of interest, they will cause no aerial congestion. Thus, we model all drones traveling along aerial path $\p\in\Paths^A$ to have a latency $\latency_\p^A$ equal to a predefined constant. This gives us the following expression for the sum of latencies experienced by parcels delivered via drones:
\begin{equation}
    \Latency^A = \sum_{\p\in\Paths^A}\flow_\p^D\latency_\p^A,
    \label{eq:drone_latency}
\end{equation}
where we no longer need to sum over the edges in a given path as done in \eqref{eq:truck_latency} since each aerial path consists of only one edge.\par

It is important to note that we are not considering the possibility of failed deliveries among the road and aerial networks. This extra consideration can be implemented into our formulation using a probabilistic approach, and defining a constant $\mu_\link$ to denote the probability of failure along a given link $\link\in\Links^R$ within the road network. The same can be done for the aerial paths $\p\in\Paths^A$.\par 
\subsection{Societal Cost}
Due to their recurrent stops, utilizing a large amount of delivery trucks can negatively impact traffic congestion as a whole, especially when considering package deliveries in dense urban environments \cite{pavone_24,impact_truck1,impact_truck2,air_pollution,parking}. To capture this effect, we define the societal cost $\societalcost$, as the normalized latency experienced by the nominal flow of traffic in units of hours: 
\begin{equation}
    \societalcost = \frac{1}{\beta}\sum\limits_{\link\in\Links^R}\flow_\link^0\latency_\link^R,
    \label{eq:societalcost}
\end{equation}
where $\beta$ represents the total flow of nominal vehicles within the network in units of vehicles per hour.
Since the edge latencies $\latency_\link$ inside the summation of \eqref{eq:societalcost} are weighted by the corresponding nominal flows $\flow_\link^0$, we are giving higher weight to roads that have a larger amount of affected drivers. We normalize this sum by the total nominal flow $\beta$ so that our societal cost $\societalcost$ represents a normalized latency per driver rather than an aggregate sum. We use the expressions "societal cost" and "societal latency" interchangeably when referring to $\societalcost$.\par 
\subsection{Operational Cost}
Finally, we model the cost of operating this delivery system, which we denote with $\Cost$ in units of dollars per hour. We break up the cost into two components associated with the road and aerial network. To calculate these components, we define $\cost^T$ and $\cost^D$ as the hourly cost of operating one truck and drone respectively. Because $\tfrac{1}{\ppertruck}\demand_\node^R$ and $\demand_\node^A$ represent node $\node$'s total hourly demand of trucks and drones, we can multiply these vehicle demands with their respective hourly rates to derive the operational cost:
\begin{equation}
    \Cost = \tfrac{\cost^T}{\ppertruck}\sum\limits_{\node\in\Nodes}\demand_\node^R + \cost^D\sum\limits_{\node\in\Nodes}\demand_\node^A.
    \label{eq:operational_cost}
\end{equation}\par
\begin{figure*}[!b]
\vspace{-10pt}
  \centering
  \subfigure[]{\includegraphics[scale=0.25]{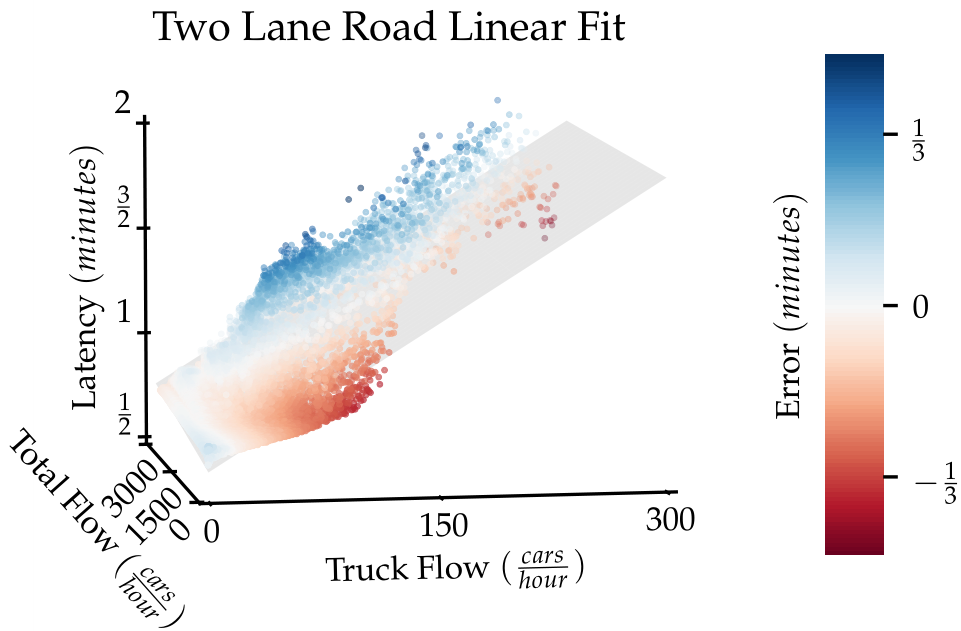}}\qquad
  \subfigure[]{\includegraphics[scale=0.25]{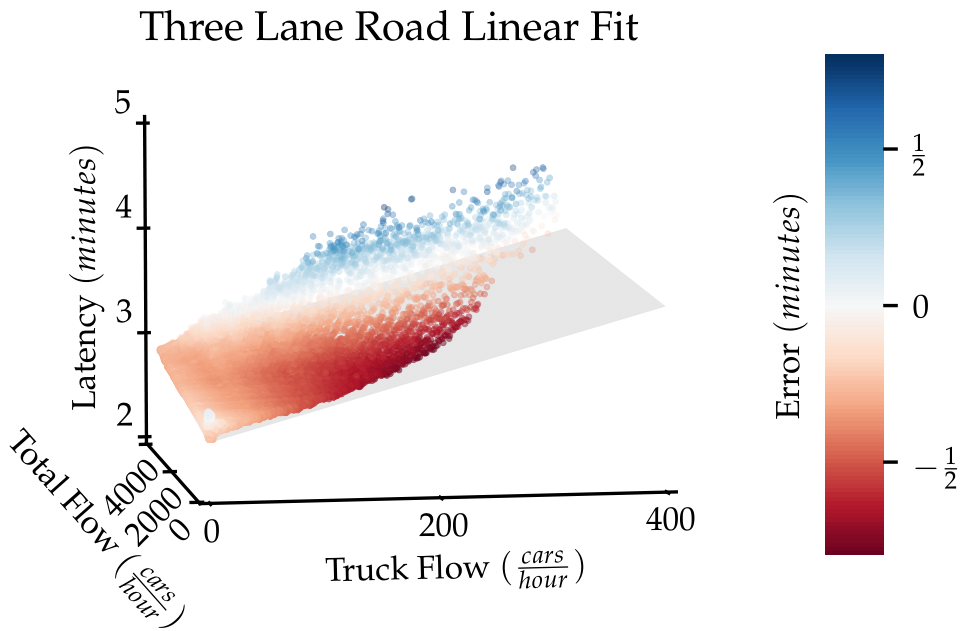}}
  \caption{The two plots (a) and (b) depict the best fit linear planes derived for our two-lane road and three-lane road respectively. The z-dimension represents latency, while the xy-dimensions correspond to our two flow components. Colors are used to visualize the error between the simulation results and the best fit planes. We can see that our affine transformations slightly over-estimate the latency when total flow is small, and under-estimate the latency when total flow is large. We emphasize that, because trucks account for a small portion of the total flow, we operate in the regime where errors are small.}
  \label{fig:sumo_3d}
\end{figure*}
\section{Optimization} \label{sec:Optimization}
Since we are working with a routing problem, our main goal is to find the optimal allocation of flow for our network given an objective function to minimize. Our decision variables for the optimization are the flows of trucks $\flow_\p^T$ along the given paths $\p\in\Paths^R$. As mentioned before, the flows of drones $\flow_\p^D$ along the aerial paths $\p\in\Paths^A$ can be seen as linear functions of our decision variable and hence do not have to be explicitly controlled in the optimization. We will now explain the objective function used, as well as the constraints required by our model.\par

\subsection{Objective Function}
We want to find the optimal routing for a parcel delivery system that minimizes both parcel latency and societal costs. Because of this, our objective function consists of two components, the average latency experienced by parcels $\Latency$, and the societal cost incurred on the system by the parcel delivery service $\societalcost$. We express our objective function as:
\begin{equation}
    J(\flow^T) = \gamma \Latency + (1-\gamma)\societalcost,
    \label{eq:obj_fun}
\end{equation}
where we use $0<\gamma<1$ as a weight that determines the relative importance of the parcel latency and societal cost objectives.\par 

\subsection{Constraints}
The first constraint of our model is the operational cost defined in \eqref{eq:operational_cost}. This is seen as an inequality constraint: $\Cost(\flow^T)\leq\Cost_0$, where $\Cost_0$ is a constant denoting the maximal operational cost allowed in dollars per hour. Additionally, we need to fulfill the system's parcel demand while maintaining conservation of flow within the network by satisfying \eqref{eq:vertex_demand} and \eqref{eq:flow_cons}. 
\subsection{Overall Optimization Problem}
Our overall optimization problem can be written as:
\begin{equation}
\begin{split}
    \min_{\flow^T}\;  J(\flow^T) &= \gamma \Latency + (1-\gamma)\societalcost \\\
    \textrm{subject to }\; \flow_\p^T&\geq0 \qquad \forall\p\in\Paths^R \:,\\
    \Cost(\flow^T) &\leq \Cost_0 \:,\\
    \demand_\node^R(\flow^T) &= \ppertruck\Big(\sum_{\link\in\uplink_\node^R}\flow_\link^T - \sum_{\link\in\downlink_\node^R}\flow_\link^T\Big) \:,\\
    \textrm{and }\; \demand_\node^A(\flow^T) &=  \demand_\node - \demand_\node^R(\flow^T), \qquad \forall\node\in\Nodes\:. 
\end{split}\label{eq:final_optimization}
\end{equation}
We will show in Section \ref{sec:Results} how the above formulation in \eqref{eq:final_optimization} can be expressed as a quadratic convex objective function with linear equality and inequality constraints.\par
\section{SUMO Simulations} \label{sec:SUMO}
Using the open source traffic simulation package ``SUMO" \cite{SUMO2018}, we derive expressions for latency as a function of the stopping and non-stopping flow components on a given link. In this section, we will briefly describe the simulation setup used and subsequently show the analyzed results.\par 
\subsection{Setup}
For all of our simulations, we used two vehicle types - cars that do not stop at delivery locations, and trucks that stop at delivery locations. The vehicle parameters used throughout the simulations are listed in Table \ref{table: veh_params}.\par
\begin{table}[h]
\caption{Vehicle Parameters}
\vspace{-5pt}
\label{table: veh_params}
\begin{center}
\begin{tabular}{|c||c||c|}
\hline
Vehicle & Car & Truck\\
\hline
Length $m$ \cite{truck_len}& $4.00$ & $7.82$\\
\hline
Max. Speed $\tfrac{km}{h}$ \cite{BOKARE20174733}& $160.93$ & $96.56$\\
\hline
Max. Acceleration $\frac{m}{s^2}$ \cite{BOKARE20174733}& $2.87$ & $1.00$\\
\hline
Max. Deceleration $\frac{m}{s^2}$ \cite{BOKARE20174733}& $4.33$ & $0.88$\\
\hline
\end{tabular}
\end{center}
\vspace{-15pt}
\end{table}
\begin{figure}[!t]
\vspace{5pt}
  \centering
  \includegraphics[scale=1]{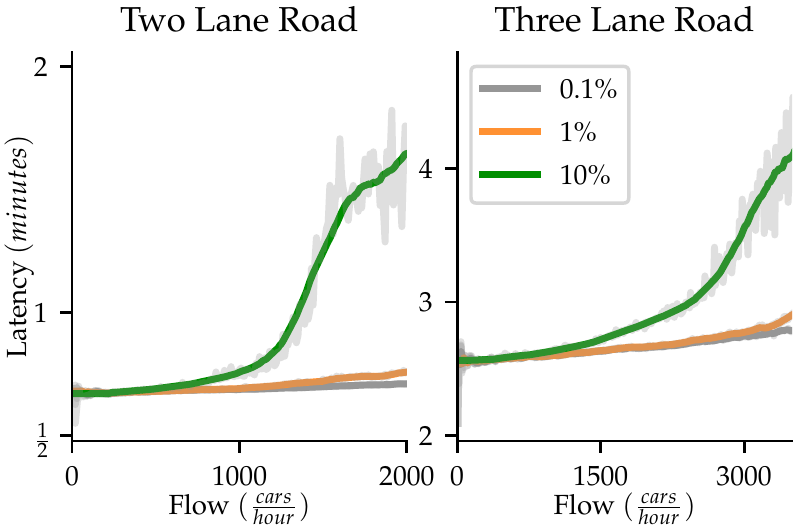}
  \caption{The three lines in each plot represent simulation results for different ratios of truck flow to total flow. As the total flow of vehicles increases, so does the latency, with larger ratios contributing to sharper increases. From the two subplots we can see that two lane roads are proportionally more affected by stopping trucks as compared to three lane roads.}
  \label{fig:sumo_latency}
  \vspace{-10pt}
\end{figure}
With these vehicle types we simulated two variants of roads, one with two lanes and the other with three lanes, placing equally spaced truck stops along both roads. For each road, we simulated over a set of two parameters - the ratio of stopping to non stopping vehicles and the total flow of vehicles. We varied the ratio from $0.001$ to $0.1$ since that is the regime of practical interest, and for each value we increased the total flow in increments of $10$ cars per hour until reaching a prescribed maximum value. The road parameters used throughout the simulations are listed in Table \ref{table: road_params}. After collecting the data for the following simulations, we analyzed our results to see how latency was affected by the input parameters.\par
\begin{table}[h!]
\caption{Road Parameters}
\vspace{-5pt}
\label{table: road_params}
\begin{center}
\begin{tabular}{|c||c||c|}
\hline
Road & Two Lane & Three Lane\\
\hline
Length $m$ & $500$ & $2000$\\
\hline
Speed Limit $\tfrac{km}{h}$ & $48.28$ & $80.47$\\
\hline
Number of Stops & $25$ & $100$\\
\hline
Duration $minutes$ & $60$ & $90$\\
\hline
Max. Flow $\tfrac{vehicles}{hour}$ & $3000$ & $5000$\\
\hline
\end{tabular}
\end{center}
\vspace{-15pt}
\end{table}
\subsection{Analysis}
For each simulation, we looked at the travel times experienced by all vehicles, averaging this value to compute the road latency given the current ratio and total flow. From Fig. \ref{fig:sumo_latency}, we can see that latency is much larger when there is a bigger portion of stopping trucks on the road, and that this effect is more severe for two lane roads. We also see that for smaller ratio values, the relationship between latency and truck flow can be modeled as linear, since in this regime each truck creates a constant amount of delay. Since in practice, we expect the flow of stopping trucks $\flow^T_\link$ along a given edge $\link$ to be much smaller than the nominal vehicular flow $\flow^0_\link$ along that edge, we consider this operating regime and approximate our edge latencies $\latency_\link^R$ as linear with respect to the two components of flow:
\begin{equation}
    \latency_\link^R = \langle \w_\link,[1,\flow_\link^R,\flow_\link^R+\flow_\link^0] \rangle,
    \label{eq:sumo_latency}
\end{equation}
where $\w_\link$ corresponds to the $3$-dimensional weight vector associated with edge $\link\in\Links^R$ and $\langle a,b\rangle$ is used to denote the inner product between vectors $a$ and $b$. As shown in Fig. \ref{fig:sumo_3d}, this is equivalent to a $3$-dimensional linear plane for each simulated road. It is important to note that our weight vectors have strictly positive terms, since latency is monotonically increasing with respect to the flow components. Using our SUMO simulations, we are able to define affine transformations that map from our network flow components to latencies.\par

\section{Simulation Studies} \label{sec:Results}
Now that we have derived the functions used to compute our edge latencies $\latency_\link^R$, we can go back to solving our overall optimization problem shown in \eqref{eq:final_optimization}. Using our expressions for average latency in \eqref{eq:avg_latency}, \eqref{eq:truck_latency}, \eqref{eq:drone_latency}, \eqref{eq:societalcost}, and our derived latency model in \eqref{eq:sumo_latency}, we can express our objective function in standard quadratic form:
\begin{equation}
    J(\flow^T,\gamma) = (\flow^T)^\intercal Q\flow^T + a^\intercal\flow^T,
    \label{eq:objective_quad}
\end{equation}
where we use $^\intercal$ to denote the transpose operation, and $\gamma$ has been absorbed into the expressions for matrix $Q$ and vector $a$.\par 
\begin{figure*}[!b]
  \centering
  \subfigure[]{\includegraphics[scale=1]{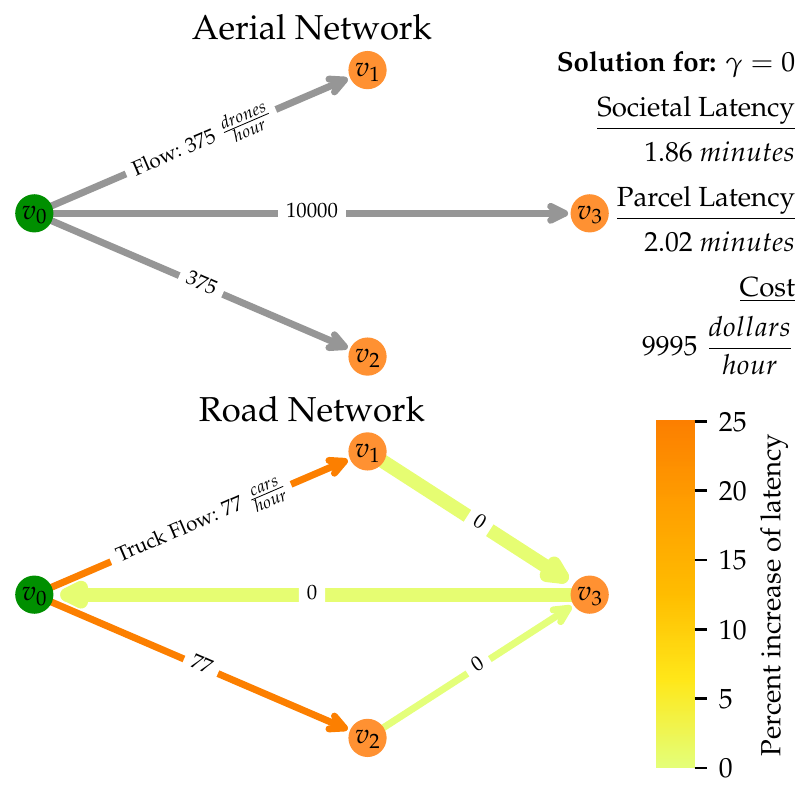}}\qquad
  \subfigure[]{\includegraphics[scale=1]{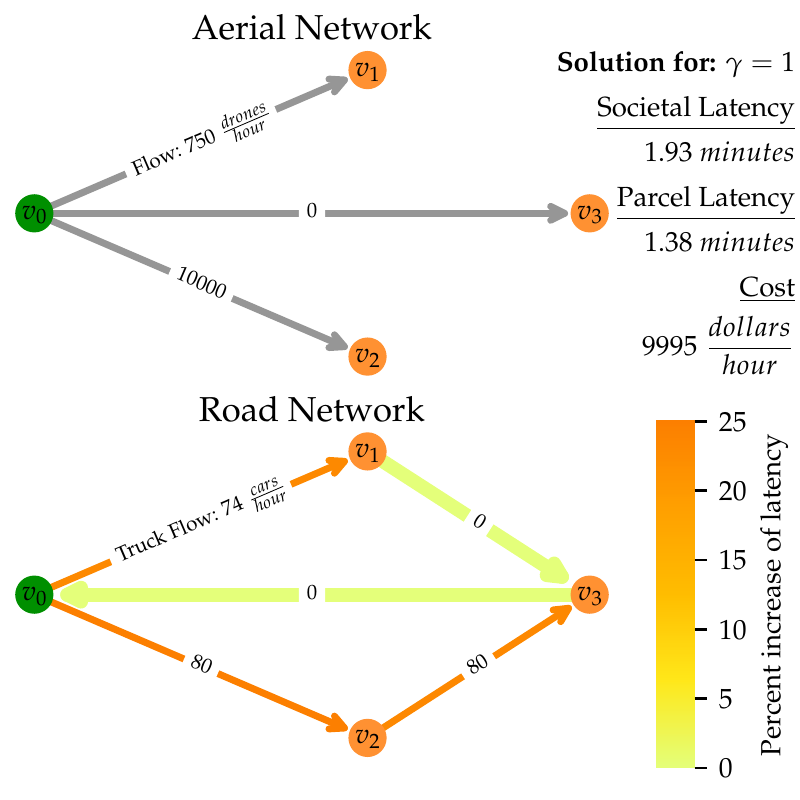}}
  \caption{The two figures (a) and (b) show the solutions to our optimization problem for different values of $\gamma$. Societal and parcel latency are calculated according to \eqref{eq:societalcost} and \eqref{eq:avg_latency} respectively. The edge labels in the aerial networks are the flows of drones along the paths. The edge labels in the road networks are the flows of trucks along the edges, which are calculated post hoc using the solution's path flows.}
  \label{fig:toy_sol}
\end{figure*} 
For a given problem setup, $Q$ is a diagonal matrix where the $i$'th row's diagonal term is $\gamma$ times the sum of the second and third components of $\w_\link$ for all $\link\in\p$, where $\p$ is the corresponding path for row $i$. Since all components of $\w_\link$ are positive, $Q$ is positive-definite. Vector $a$'s construction requires many operations, and hence is left out for readability purposes as it can be derived using the equations provided. 
We go on to show a small example that visualizes how we construct our network model and optimize. Following this, we preform a case study on the Sioux Falls transportation network, shown in \ref{fig:sioux}, to demonstrate the scalability of our methods.\par 
\subsection{Complexity and Scalability}
Since the problem posed in \eqref{eq:final_optimization} is defined by a convex quadratic objective function, it can be efficiently solved using quadratic programming (QP) techniques. Inspired by Karmarkar's  Linear Programming algorithm \cite{scale_4_karmarkar1984new}, multiple interior-point algorithms have adapted his method for the QPs, showing convergence to the optimal solution in $O(n^3L)$ arithmetic operations \cite{scale_1_goldfarb1990n,scale_2_ye1989extension,scale_3_monteiro1989interior}, where $n$ is the dimension of the decision variable and $L$ is the number of binary symbols needed to record the input information.\par
In our specific case, $n$ is the number of paths considered in the problem $|\Paths^R|$, and $L$ is the number of bits required to store the information in $Q,a$ from \eqref{eq:objective_quad} as  well as the constraints defined in \eqref{eq:final_optimization}. Since the number of paths between two nodes exhibits factorial scaling with the number of links between them, the total number of paths can scale worse than exponentially. Fortunately, we can restrict ourselves to consider paths up to a certain length, largely limiting the computational overhead. 
\subsection{Example Network}
\begin{figure}[!t]
\vspace{5pt}
  \centering
  \includegraphics[scale=1]{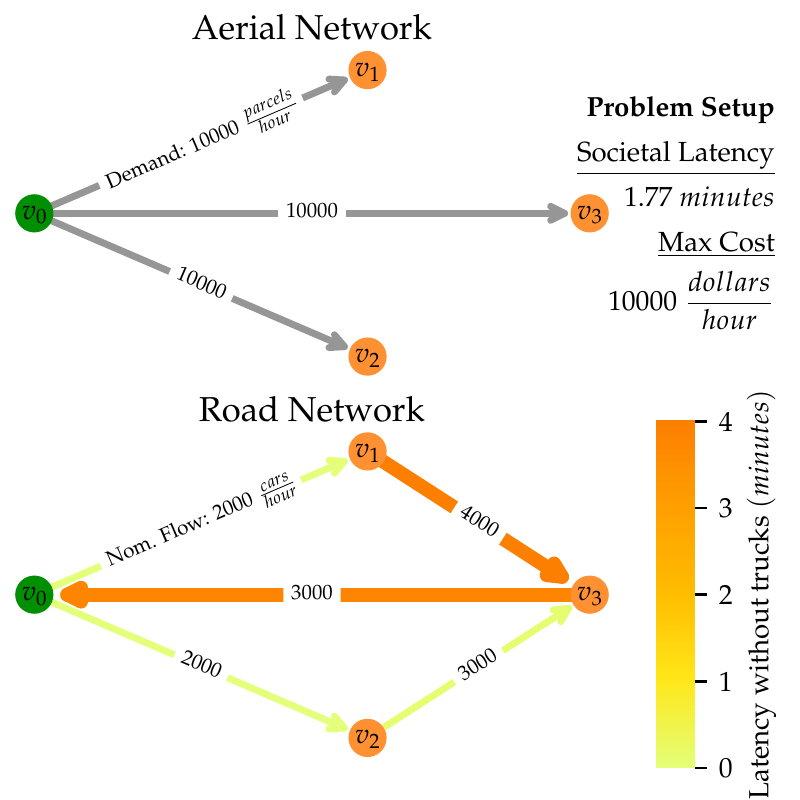}
  \caption{This figure depicts the setup for our small network example, note that it is not drawn to scale. The top graph represents $\Graph^A$, labeling the parcel demands for each node along the corresponding edge. The bottom graph represents $\Graph^R$, labeling the nominal flows of vehicles along the edges. The thicker edges represent three lane roads that are $2$ kilometers long, while the thinner edges represent two lane roads that are $0.5$ kilometers long.}
  \label{fig:toy_setup}
  \vspace{-20pt}
\end{figure}

For our example network we choose a set of four nodes $\Nodes$ and two sets of edges, $\Links^R$ and $\Links^A$, building our road and aerial networks, $\Graph^R$ and $\Graph^A$, as shown in Fig. \ref{fig:toy_setup}. The nominal flows $\flow^0_\link$ and parcel demands $\demand_\node$ are shown in Fig. \ref{fig:toy_setup}, where $\beta=14,000$ vehicles per hour is used to compute the societal latency. We use $\cost^T=30$ dollars per hour and $\cost^D=0.5$ dollars per hour as our operational costs for trucks and drones respectively, allowing each truck to carry a total of $\ppertruck=125$ parcels \cite{drone_cost}. Finally, we set the drone speed as $25$ kilometers per hour, as this should allow for the longest flight distance \cite{pavone_14}.\par
With this setup, we optimize for our objective function using multiple values of $\gamma$ to shift priority between minimizing societal latency and minimizing parcel latency. For this example, we chose to use the full set of available paths as our decision variable, depicting our results in Table \ref{table: toy_study}.\par
\begin{table}[h]
\vspace{5pt}
\caption{Example Network Results}
\vspace{-5pt}
\label{table: toy_study}
\begin{center}
\begin{tabular}{|c||c||c|}
\hline
Setup & Without Drones & With Drones\\
\hline
$\gamma=0$ & $1.98$, $2.29$ & $1.86$, $2.02$\\
\hline
$\gamma=\tfrac{1}{2}$ & $1.98$, $1.81$ & $1.89$, $1.70$\\
\hline
$\gamma=1$ & $1.98$, $1.63$ & $1.93$, $1.38$\\
\hline
\end{tabular}
\end{center}
\caption*{TABLE III: This table shows the results for our example network using multiple values of $\gamma$. Each entry includes two values: societal latency and parcel latency, both of which are reported in minutes.}
\vspace{-35pt}
\end{table} 

As expected, we can see that utilizing drones helps with minimizing both societal latency and parcel latency. More specifically, we see that with drones we can prioritize minimizing for societal latency by decreasing $\gamma$, meanwhile without drones, decreasing $\gamma$ has little effect on societal latency. This is due to the chosen setup of our example, showcasing that trucks will in general have a negative impact on congestion within the road network.\par 

The solutions for $\gamma=0$ and $\gamma=1$ are depicted in Fig. \ref{fig:toy_sol}, where the edge labels correspond to the drone and truck flows derived from our solution. We can see that setting $\gamma$ to $0$ routes the drones to $\node_3$, the farthest node. This causes an increase in parcel latency, but keeps trucks off of the road connecting $\node_2$ to $\node_3$. On the other hand, setting $\gamma$ to $1$ routes the drones to $\node_2$ and lets the trucks make faster deliveries to $\node_3$ instead. This reduces the overall parcel latency, but increases societal latency due to the added congestion on the road connecting $\node_2$ to $\node_3$.\par

\subsection{Sioux Falls Case Study}
\begin{figure}[t]
\vspace{5pt}
  \centering
  \includegraphics[width=2.5in]{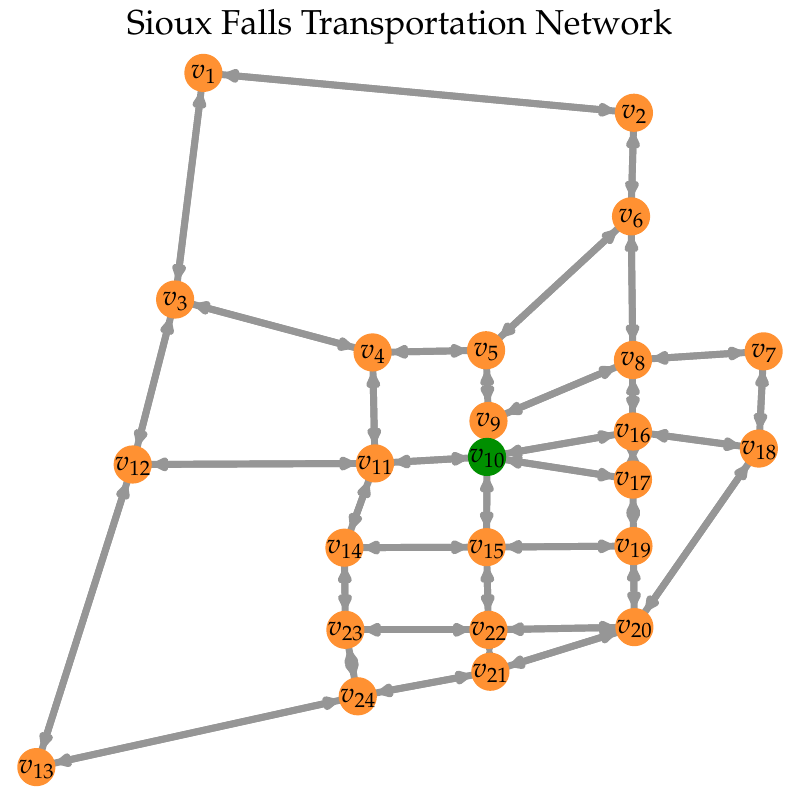}
  \caption{The Sioux Falls transportation network is depicted above. The green node $\node_{10}$ represents our delivery hub's location.}
  \label{fig:sioux}
  \vspace{-10pt}
\end{figure}

For our case study, we use the Sioux Fall transportation network shown in Fig.~\ref{fig:sioux} since it has been widely used as a standard test for traffic assignment \cite{sioux_falls}. The network contains $24$ nodes and $76$ directed edges, where we choose node $\node_{10}$ as our central delivery hub. Using the reported values for the best known edge flows and the edge capacities, we scale our nominal flows $\flow_\link^0$ to be proportional to the ratios between edge flows and edge capacities. We choose $24$ of the $76$ edges to be three lane roads since their capacities are over $10,000$ vehicles per hour, and scale the weights $\w_\link$ for all edges $\link\in\Links^R$ based on the physical distances between the nodes each edge connects. We consider all $1133$ paths that consist of eight or less edges, and optimize for our objective function using parcel demands of $5000$ parcels per hour for each node, setting the total nominal flow $\beta$ to $81,000$ vehicles per hour, and allowing a maximum operational cost of $50,000$ dollars per hour. 
We summarize the results of our case study in Table \ref{table: case_study}.\par
First of all, we see a similar trend as depicted in our example network: utilizing drones minimizes both societal latency and parcel latency. In addition, we see that drones have a greater affect on parcel latency than was depicted by our example network. Since our case study's road network consists of an actual grid rather than direct links between nodes, we see the advantages drones have due to their ability to utilize direct paths between nodes. Overall, our framework quantifies the potential drones have on mitigating road congestion and solves for the path routing needed to minimize a chosen combination of societal latency and parcel latency, while remaining computationally feasible to scale due to its reliance on convex quadratic optimization.\par 
\begin{table}[h]
\caption{Case Study Results}
\vspace{-5pt}
\label{table: case_study}
\begin{center}
\begin{tabular}{|c||c||c|}
\hline
Setup & Without Drones & With Drones\\
\hline
$\gamma=0$ & $9.33$, $10.97$ & $8.70$, $7.79$\\
\hline
$\gamma=\tfrac{1}{2}$ & $9.39$, $8.94$ & $8.82$, $7.33$\\
\hline
$\gamma=1$ & $9.51$, $8.29$ & $9.00$, $7.20$\\
\hline
\end{tabular}
\end{center}
\caption*{TABLE IV: This table shows the results for our case study using multiple values of $\gamma$. Each entry includes two values: societal latency and parcel latency, both of which are reported in minutes.}
\vspace{-25pt}
\end{table}

\section{Conclusion}
We proposed a network model for a congestion-aware Bi-modal delivery system incorporating drones. Using our simulation results, we constructed a convex optimization framework which can solve for the optimal routing of our delivery system. Finally, we showed the scalability of our framework by applying it on the well studied Sioux Falls transportation network. Below we discuss the limitations of our work, including directions for possible improvement.\par

An inaccuracy of our framework is the planar fit latency models depicted in Fig. \ref{fig:sumo_3d}. One way to address this issue is by utilizing a piece-wise linear fit instead. This addition would make our main objective function piece-wise quadratic, meaning we can no longer rely on simple quadratic programming solvers for our optimization. We point interested readers to a study of such piece-wise quadratic solvers, leaving this extension to future work \cite{pq_study}.\par 
In addition to the inaccuracies of our latency estimates, one may view the simplicity of our proposed networks as a limitation as well. We discourage this perspective, since the inherent simplicity is what allows for the construction of a scalable and convex optimization framework. As pointed out before, prior works have shown lack of applicability in large scale networks due to their reliance on heuristic techniques and MILP formulations to solve for routing policies \cite{pavone_8,pavone_37_uav_route,drone_vrp}. On the other hand, our proposed network model and optimization framework lets us scale to arbitrarily large networks, allowing delivery systems to find routing policies for their trucks and drones, and evaluate the effects different routing policies can have on both their parcel delivery latency and the traffic conditions of the roads.\par 
\section*{Acknowledgments}
We acknowledge funding by NSF ECCS Grant \#$1952920$.
\bibliographystyle{IEEEtran}
\bibliography{refs}
\end{document}

%% file: macros.tex
\usepackage{bm}
\usepackage[font=small]{caption}

\usepackage{titlesec}
\titlespacing*{\subsection}{0pt}{0.15\baselineskip}{0.05\baselineskip}
\titlespacing*{\section}{0pt}{0.6\baselineskip}{0.45\baselineskip}

\usepackage[colorlinks,citecolor=red]{hyperref}
\usepackage{color}

\newcommand{\Graph}{G}
\newcommand{\Nodes}{\mathcal{N}}
\newcommand{\cnode}{v_0}
\newcommand{\Links}{\mathcal{E}}
\newcommand{\node}{v}
\newcommand{\link}{e}
\newcommand{\Paths}{\mathcal{P}}
\newcommand{\p}{p}
\newcommand{\uplink}{\mathcal{U}}
\newcommand{\downlink}{\mathcal{D}}
\newcommand{\demand}{d}
\newcommand{\flow}{f}
\newcommand{\ppertruck}{m}
\newcommand{\Latency}{L}
\newcommand{\latency}{\ell}

\newcommand{\societalcost}{\Latency^S}
\newcommand{\cost}{c}
\newcommand{\Cost}{C}
\newcommand{\w}{\omega}

